\documentclass[11pt]{amsart}
\usepackage{amsmath}
\usepackage{amssymb, verbatim}
\usepackage{pstricks,pst-node,pst-plot}
\usepackage{comment}
\usepackage{color}
\usepackage{extarrows}
\usepackage{graphicx}

\usepackage[linktocpage=true,hypertexnames=false]{hyperref}
\hypersetup{colorlinks=true, citecolor=blue, linkcolor=blue, urlcolor=blue, pdfauthor={}, pdftitle={},breaklinks=true}

\usepackage[numbers]{natbib}

\title[]{Calabi-Yau threefolds across quadratic singularities}

\theoremstyle{plain}
\newtheorem{thm}{Theorem}[section]

\newtheorem{defn}[thm]{Definition}

\newtheorem{cor}[thm]{Corollary}
\newtheorem{conj}[thm]{Conjecture}

\theoremstyle{definition}
\newtheorem{ex}[thm]{Example}
\newcommand{\exampleend}{\hfill $\blacktriangle$}
\newtheorem{rk}[thm]{Remark}

\numberwithin{equation}{section}

\newcommand{\be}{\begin{equation}}
\newcommand{\bea}{\begin{eqnarray}}
\newcommand{\eea}{\end{eqnarray}} 
\newcommand{\ee}{\end{equation}}

\renewcommand{\leq}{\leqslant}
\renewcommand{\geq}{\geqslant}

\renewcommand{\epsilon}{\varepsilon}
\renewcommand{\phi}{\varphi}

\begin{document}

\author[S. Picard]{S\'ebastien Picard}
  \email{spicard@math.ubc.ca}
  \address{Department of Mathematics, The University of British
    Columbia, 1984 Mathematics Road, Vancouver BC Canada V6T 1Z2}
\thanks{The author is supported by an NSERC Discovery Grant.}

\maketitle

\begin{abstract}
These are lecture notes on non-K\"ahler complex threefolds presented at the MATRIX program ``The geometry of moduli spaces in string theory''. We review some basics of Calabi-Yau geometry in Section 1, describe topological features of the conifold transition in Section 2, and survey recent developments on the geometrization of conifold transitions in Section 3.
\end{abstract}


  \section{Calabi-Yau Threefolds}

  
  \subsection{Definitions}
 The topic of this survey is the Calabi-Yau threefold.
  
  \begin{defn}
We take the definition of a Calabi-Yau threefold $X$ to be a compact complex manifold with $\dim_{\mathbb{C}} X = 3$ which is projective, satisfies $H^1(X,\mathbb{C})=0$, and admits a holomorphic volume form $\Omega$.
    \end{defn}

    Recall that a holomorphic volume form on $X$ is given by a $(3,0)$-form $\Omega$ such that
    \[
\Omega \overset{\text{loc}}{=} f(z) \, dz^1 \wedge d z^2 \wedge dz^3
    \]
    where $f(z)$ is a local nowhere vanishing holomorphic function.

    \begin{rk}
The definition of a Calabi-Yau threefold is not standardized in the literature, and some setups do not require the vanishing of $H^1(X,\mathbb{C})$ or generalize the existence of a holomorphic volume form $\Omega$ to $c_1(X)=0$. A consequence of the Calabi-Yau \cite{Calabi, Yau78} theorem is that a K\"ahler manifold $X$ with $c_1(X)=0$ must have $K_X$ holomorphically torsion, but not necessarily trivial. Also, in the given definition projectivity is redundant as it follows from $H^1(X,\mathbb{C})=0$, which implies $h^{0,2}=0$, and projectivity then follows from the Kodaira embedding theorem.
      \end{rk}

      \begin{ex}
        A basic example of a Calabi-Yau threefold is the following quintic threefold:
        \[
X = \bigg\{ \sum_{i=0}^4 Z_i^5 = 0 \bigg\} \subset \mathbb{P}^4.
        \]
This defines a smooth complex manifold by the implicit function theorem, since a hypersurface $\{P=0\} \subseteq \mathbb{P}^n$ with $P$ a homogeneous polynomial is smooth if there is no non-zero point where simultaneously $P=0$ and $DP = 0$. The holomorphic volume form in this case is
        \[
\Omega \overset{\text{loc}}{=} \frac{dw_1 \wedge d w_2 \wedge d w_3}{w_4^4}, \quad w_i = \frac{Z_i}{Z_0}
\]
over the open set $U= \{ Z_0 \neq 0, \ w_4 \neq 0 \}$. This formula can be verified to glue on overlaps of similar local open sets to define a global section $\Omega \in H^0(X,\Omega^3)$. That $X$ satisfies $H^1(X,\mathbb{C})=0$ follows from the Lefschetz hyperplane theorem.  \exampleend
        
\end{ex}

\begin{ex}
  Any homogeneous degree 5 polynomial $P$ will also define a Calabi-Yau threefold $\{P = 0 \} \subseteq \mathbb{P}^4$ provided the zero locus is smooth.
  \exampleend
\end{ex}

\begin{ex}
For more examples of Calabi-Yau threefolds beyond quintics $\{P = 0 \} \subseteq \mathbb{P}^4$, we refer to for example H\"ubsch's book \cite{Hubsch}. \exampleend 
  \end{ex}

      The Hodge diamond of a Calabi-Yau threefold has two parameters:

      \[
\begin{array}{ccccccc}
  &   &   & 1 &   &   &   \\
  &   & 0 &   & 0 &   &   \\
  & 0 &   & h^{1,1} &   & 0 &   \\
1 &   & h^{2,1} &   & h^{2,1} &   & 1 \\
  & 0 &   & h^{1,1} &   & 0 &   \\
  &   & 0 &   & 0 &   &   \\
  &   &   & 1 &   &   &   
\end{array}
\]
Mirror symmetry \cite{Dixon, GreenePlesser, COGP, LVW, Kontsevich95, SYZ} predicts that Calabi-Yau threefolds come in pairs $(X,\check{X})$ exchanging the two parameters $(h^{1,1},h^{2,1})$.
\[
h^{1,1}(X) = h^{2,1}(\check{X}), \quad h^{1,1}(\check{X}) = h^{2,1}(X).
\]
A first indication of the role of non-K\"ahler manifolds comes from the inherent asymmetry of Calabi-Yau threefolds where $h^{2,1}$ may vanish while $h^{1,1}$ cannot. Indeed, a K\"ahler metric $\omega > 0$ creates a non-zero K\"ahler class $[\omega] \in H^{1,1}(X)$. Kontsevich \cite{Kontsevich} has suggested that the mirror theory of curve enumeration on a threefold $X$ with $h^{2,1}(X)=0$ should involve Hodge structures on a non-K\"ahler complex threefold $\check{X}$ with $h^{1,1}(\check{X})=0$.

We now review the significance of the two parameters $(h^{1,1},h^{2,1})$.

\subsection{Discussion of $h^{1,1}$}
Calabi-Yau threefolds are studied in differential geometry because they support solutions to the Ricci-flat equation. Recall that on a Riemannian manifold $(M,g)$, the Riemann curvature tensor is a second order invariant of $g_{ij}$ given by
\[
R_{pq}{}^k{}_j = \partial_p \Gamma_q{}^k{}_j + \Gamma_p{}^k{}_r \Gamma_q{}^r{}_j - (p \leftrightarrow q),
\]
with
\[
\Gamma_i{}^k{}_j = \frac{g^{k p}}{2} (-\partial_p g_{ij} + \partial_i g_{pj} + \partial_j g_{ip}).
\]
The Ricci tensor is given by
\[
R_{ij} = - R_{ik}{}^k{}_j,
\]
and a fundamental equation in geometry and physics is the Ricci-flat equation
\[
  R_{ij}=0.
\]
On a Calabi-Yau threefold $X$, one can look for solutions of the following form. Choose a background reference K\"ahler metric $g_{\mu \bar{\nu}}$ on $X$ such as the pullback of the Fubini-Study metric in the embedding $\iota: X \rightarrow \mathbb{P}^N$, and set the ansatz
\begin{equation} \label{Kclass}
\tilde{g}_{\mu \bar{\nu}} = g_{\mu \bar{\nu}} + \partial_\mu \partial_{\bar{\nu}} u > 0
\end{equation}
for an unknown potential function $u \in C^\infty(X)$. Let Greek indices denote holomorphic coordinates $\{ z^\alpha \}$ on the complex manifold $X$. K\"ahler \cite{Kahler} computed the Ricci tensor for ansatz \eqref{Kclass} and derived the ``sehr elegant'' equation
\[
R_{\alpha \bar{\beta}} = \partial_\alpha \partial_{\bar{\beta}} \log |\Omega|_{\tilde{g}}^2, \quad R_{\alpha \beta} = 0
\]
where
\[
  |\Omega|^2_{\tilde{g}} \overset{{\rm loc}}{=} \frac{f(z) \overline{f(z)}}{\det \tilde{g}_{\mu \bar{\nu}}}.
\]
We can then look for solutions of the Ricci-flat equation by setting $|\Omega|_{\tilde{g}} = 1$, so that the equation to solve becomes
\begin{equation} \label{CMA}
\det(g_{\alpha \bar{\beta}} + u_{\alpha \bar{\beta}}) = e^f \det g_{\alpha \bar{\beta}}
\end{equation}
where $e^f = |\Omega|^2_g$ is given. This geometric complex Monge-Amp\`ere equation is an adaptation of the fundamental PDE
\[
\det D^2 u = \psi >0,
\]
on a domain in $\mathbb{R}^n$. The complex Monge-Amp\`ere equation \eqref{CMA} was solved by Yau.

\begin{thm} [Yau's theorem \cite{Yau78}] There exists a unique solution $u \in C^\infty(X)$ to \eqref{CMA} with $\sup_X u = 0$.
\end{thm}

Consequently, there exists Ricci-flat metrics on a Calabi-Yau threefold $X$. The ansatz \eqref{Kclass} is interpreted as finding a unique K\"ahler Ricci-flat metric in a given K\"ahler class $[\omega] \in H^{1,1}(X,\mathbb{R})$, where a K\"ahler metric $g$ defines a form $\omega \in \Omega^{1,1}(X,\mathbb{R})$ via $\omega = i g_{\mu \bar{\nu}} \, dz^\mu \wedge d \bar{z}^\nu$ and so by the $\partial \bar{\partial}$-lemma the ansatz \eqref{Kclass} is equivalent to setting $[\tilde{\omega}] = [\omega]$.

The construction comes in families parametrized by the choice of background K\"ahler metric $g_{\mu \bar{\nu}}$, but as noted above really the parameter is the K\"ahler class $[\omega] \in H^{1,1}(X,\mathbb{R})$. Thus $h^{1,1}$ counts the dimension of the space of K\"ahler Ricci-flat metrics on $X$ with fixed complex structure.

\subsection{Discussion of $h^{2,1}$}
Having discussed the parameter $h^{1,1}$, we now give the interpretation of the parameter $h^{2,1}$.

Recall that a family of complex manifolds over a base $B \subseteq \mathbb{C}^k$ with $0 \in B$ is defined as a proper holomorphic submersion
  \[
\pi: \mathcal{X} \rightarrow B
  \]
  where $\mathcal{X}$ is a complex manifold. The fibers $X_b = \pi^{-1}(b)$ for $b \in B$ are then all compact complex manifolds varying holomorphically in the parameter $b$. Ehresmann's theorem (see e.g. \cite{Kodaira} for an exposition) states that after possibly replacing $B$ with a neighborhood of $0$, there exists a diffeomorphism
  \[
\Psi:  X_0 \times B \overset{\cong}{\rightarrow} \mathcal{X}
  \]
  such that $\pi \circ \Psi = {\rm pr}_2$, with ${\rm pr}_2$ the projection to the second factor. Therefore all manifolds $X_b$ are diffeomorphic. Each $X_b$ comes with a complex structure tensor $J_b$, and via the diffeomorphism we obtain a family of complex structures also denoted $J_b$ over the fixed differentiable manifold $X_0$.

\begin{rk}
Recall that given a complex manifold $X$ with holomorphic coordinates $\{ z^\alpha \}$, the corresponding complex structure tensor $J \in \Gamma({\rm End} \, T_{\mathbb{C}} X)$ is defined by
\[
J^\alpha{}_\beta = i \delta^\alpha{}_\beta \quad J^{\bar{\alpha}}{}_{\bar{\beta}} = - i \delta^\alpha{}_\beta, \quad J^\alpha{}_{\bar{\beta}} = J^{\bar{\alpha}}{}_\beta = 0,
\]
in a holomorphic coordinate system. Here we use index notation for components of a tensor $A \in \Gamma({\rm End} \, T_{\mathbb{C}} X)$ so that e.g.
\[
A \bigg( \frac{\partial}{\partial z^\beta} \bigg) = A^\alpha{}_\beta \frac{\partial}{\partial z^\alpha} + A^{\bar{\alpha}}{}_{\beta} \frac{\partial}{\partial \bar{z}^{\bar{\alpha}}}.
\]
In a smooth real coordinate system, a complex structure tensor $J \in \Gamma({\rm End} \, TX)$ is characterized by the property
\begin{equation} \label{cplxstr}
J^2 = - id, \quad N(J) = 0
\end{equation}
where $N \in \Omega^2(TM)$ is the Nijenhuis tensor given by
\[
N^k{}_{ij} = {1 \over 4} \bigg(  J^k{}_p \partial_j J^p{}_i + \partial_p J^k{}_j J^p{}_i - (i \leftrightarrow j) \bigg).
\]
The Newlander-Nirenberg theorem (see e.g. \cite{Demailly}) states that the existence of a $J \in \Gamma({\rm End} \, TX)$ satisfying \eqref{cplxstr} on a smooth manifold $X$ is equivalent to the existence of holomorphic coordinate charts. Such a $J$ splits $T_{\mathbb{C}} X = T^{1,0}X \oplus T^{0,1}X$ into the $+i$ and $-i$ eigenvalues of $J$, and the interpretation of $N$ is that $N=0$ if and only if $[U,V] \in T^{1,0} X$ for all $U,V \in \Gamma(T^{1,0}X)$.
\end{rk}

Given a family of complex structures $J_t$ varying smoothly with a parameter $t \in (-\epsilon, \epsilon)$ on a fixed smooth manifold $X$, differentiation along the path defines the fluctuation tensor
\[
\eta = {d \over dt} \Big|_{t=0} J_t.
\]
One can verify the following calculations :

\begin{enumerate}
\item Differentiating $J_t^2 = - id$ shows that
  \[
\eta \in \Omega^{0,1}(T^{1,0}X),
  \]
  so that for example $\eta^{\bar{\alpha}}{}_{\bar{\beta}}=0$ in holomorphic coordinates at $t=0$.

\item Differentiating $N(J_t)=0$ leads to the constraint
  \[
\bar{\partial} \eta = 0,
  \]
  so that $[\eta]$ defines a class in $H^{0,1}(X,T^{1,0}X)$.
  
\item If the family $J_t$ comes from $J_t = (\Theta_t)^{-1}_* J_0 (\Theta_t)_*$ where $\Theta_t: X \rightarrow X$ is a smoothly varying family of diffeomorphisms, then
  \[
\eta = \bar{\partial} V, \quad V \in \Gamma(T^{1,0}X)
  \]
  so that $[\eta] = [0]$ in $H^{0,1}(X,T^{1,0}X)$.
\end{enumerate}

In summary, a family of complex structures produces a cohomology class in $H^{0,1}(X,T^{1,0}X)$, and deformations just coming from families of diffeomorphisms are not counted.

\begin{rk}
  On a Calabi-Yau threefold, $K_X \cong \mathcal{O}_X$ and this parameter space has dimension
\[
h^{0,1}(T^{1,0}X) = h^{0,1}(T^{1,0} X \otimes K_X) = h^{2,1}.
\]
For further analysis of the parameter space of complex structures on a Calabi-Yau threefold, see Candelas-de la Ossa \cite{CdlO91}.
\end{rk}

The natural question that arises is the inverse problem: does a given class $[\eta] \in H^{0,1}(X,T^{1,0}X)$ come from a family $(X,J(t))$ of complex structures with $[{d \over dt}|_{t=0} J ]= [ \eta] $? It does, and this is known as the BTT theorem (see e.g. \cite{PopoBook} for a recent exposition).

\begin{thm} [Bogomolov-Tian-Todorov Theorem \cite{Tian87,Todorov}] Let $X$ be a compact K\"ahler manifold admitting a holomorphic volume form.  There is a family of complex manifolds
  \[
    \pi: \mathcal{X} \rightarrow B, \quad \pi^{-1}(0)= X
  \]
  where $B$ is a neighborhood of the origin in $H^{0,1}(X,T^{1,0}X)$. Hence given any $[\eta] \in H^{0,1}(X,T^{1,0}X)$, there exists a family $J_t$ of complex structures varying smoothly with $t \in (-\epsilon,\epsilon)$ such that  $[{d \over dt} \big|_{t=0} J ]= [ \eta]$.

\end{thm}

\subsection{Rational curves}

A key topic in the study of Calabi-Yau threefolds is its set of rational curves. These special submanifolds can be studied from various different points of view. For example, the extraordinary work of Candelas-de la Ossa-Green-Parkes \cite{COGP} shocked the algebraic geometry community by predicting a formula for the number of rational curves of degree $d$ by methods of string theory.

Let $X$ be a compact complex threefold with holomorphic volume form. Let $C \subseteq X$ be a smooth complex submanifold of complex dimension 1 with $C \cong \mathbb{P}^1$. The exact sequence defining the normal bundle is
\[
0 \rightarrow TC \rightarrow TX|_C \rightarrow \mathcal{N}_C \rightarrow 0.
\]
The first Chern class of this sequence satisfies
\[
c_1(TX|_C) = c_1(TC) + c_1(\mathcal{N}_C)
\]
and so $c_1(\mathcal{N}_C) = -2$. Grothendieck's classification of holomorphic bundles over $\mathbb{P}^1$ implies that the normal bundle of $C$ must be isomorphic to
\[
\mathcal{O}(a) \oplus \mathcal{O}(-a-2) \rightarrow \mathbb{P}^1.
\]
The case of $\mathcal{O}(-1)^{\oplus 2}$ is distinguished by its symmetry and that it is the only case in which the curve is rigid, as the other bundles admit infinitesimal deformations in $H^0(\mathcal{N}_C)$.

Such $(-1,-1)$-curves are expected to exist in abundance on Calabi-Yau threefolds \cite{Friedman,Kontsevich}. Work of Clemens \cite{ClemensIHES} and Katz \cite{Katz} guarantees the existence of $(-1,-1)$-curves on the generic quintic threefold. 

  \begin{ex} Here is an explicit example of a $(-1,-1)$-curve following \cite{Katz}. Consider the quintic threefold given by
    \[
X = \{ f_2 Z_2 + f_3 Z_3 + f_4 Z_4 = 0 \} \subseteq \mathbb{P}^4,
    \]
and the embedded holomorphic curve
    \[
      C = \{ Z_2=Z_3=Z_4=0 \} \subseteq X,
    \]
 where we take the specific $f_i$ to be
    \[
X = \bigg\{ (Z_0^4 + Z_2^4) Z_2 + (Z_0^2 Z_1^2 + Z_3^4) Z_3 + (Z_1^4 + Z_4^4) Z_4 = 0 \bigg\} \subseteq \mathbb{P}^4,
\]
so that $X$ is smooth. We will compute the normal bundle of $C \subset X$. Recall the definition of the normal bundle: let $p \in C$, and suppose $p \in U \subset X$ is a coordinate chart such that
    \[
      C \cap U = \{y = 0 \}
    \]
    where $(x_1,y_1,y_2)$ are holomorphic coordinates on $U$. Suppose also $p \in \tilde{U}$ with another such submanifold coordinate system $(\tilde{x}_1, \tilde{y}_1,\tilde{y}_2)$. Then on overlaps $U \cap \tilde{U} \cap C$, the transition function
    \[
 {\partial \tilde{y} \over \partial y} \bigg|_C
    \]
    defines the data of the normal bundle $\mathcal{N}_C$.
    \begin{itemize}

      \item Let $U$ be an open set in $\{ Z_0 \neq 0 \} \subseteq X$ near $C$. Here we use local coordinates on $X$ given by $(w_1,w_3,w_4)$ where $w_i = Z_i/Z_0$. The equation of $X$ is
    \begin{equation} \label{implicitdiff}
(1+w_2^4)w_2 + (w_1^2+w_3^4)w_3 + (w_1^4+w_4^4)w_4=0,
\end{equation}
and the curve appears as
    \[
C \cap U = \begin{cases} 
w_3 = 0 \\
w_4 = 0.
\end{cases}
\]
By the inverse function theorem, $w_2(w_1,w_3,w_4)$ is a holomorphic function of the remaining coordinates.

   \item Let $\tilde{U}$ be an open set in $\{ Z_1 \neq 0 \} \subseteq X$ near $C$. Here we use local coordinates on $X$ given by $(\tilde{w}_0,\tilde{w}_2,\tilde{w}_3)$ where $\tilde{w}_i = Z_i/Z_1$. The curve appears as
    \[
C \cap \tilde{U} = \begin{cases} 
\tilde{w}_2 = 0 \\
\tilde{w}_3 = 0.
\end{cases}
\]

\item We have covered the curve by open charts: $C \subseteq U \cup \tilde{U}$. We conform to earlier conventions by setting $(x,y_1,y_2)=(w_1,w_3,w_4)$ and $(\tilde{x}, \tilde{y}_1, \tilde{y}_2)=(\tilde{w}_0, \tilde{w}_2, \tilde{w}_3)$. Then we recognize $x$ and $\tilde{x}$ as the coordinates along $C$ and since $\tilde{x} = x^{-1}$ we see $C = \mathbb{P}^1$. The transition function of the normal bundle may be computed by implicit differentiation of \eqref{implicitdiff} which gives
    \[
\frac{\partial \tilde{y}}{\partial y} \bigg|_C = \begin{bmatrix} -x & -x^3 \\ x^{-1} & 0 \end{bmatrix}.
    \]
    This transition function is a disguise of $\mathcal{O}(-1)\oplus \mathcal{O}(-1)$. Recall that two sets of transition functions $\{ g_{\alpha \beta}, U_\alpha \cap U_\beta \}$ and $\{ \hat{g}_{\alpha \beta}, U_\alpha \cap U_\beta \}$ define isomorphic bundles if there exists $\lambda_\alpha : U_\alpha \rightarrow GL(k,\mathbb{R})$ with
    \[
\hat{g}_{\alpha \beta} = \lambda_\alpha g_{\alpha \beta} \lambda_\beta^{-1} .
    \]
    In this example, one can find $2 \times 2$ matrices $\lambda_1$, $\lambda_2$ such that
    \[
\begin{bmatrix} x & 0 \\ 0 & x \end{bmatrix} = \lambda_1 \begin{bmatrix} -x & -x^3 \\ x^{-1} & 0 \end{bmatrix} \lambda_2^{-1}
    \]
    with the matrix on the left the familiar transition function defining $\mathcal{O}(-1)\oplus \mathcal{O}(-1) \rightarrow \mathbb{P}^1$.\exampleend
    \end{itemize}
  \end{ex}

  We will later use a local model for the holomorphic curve $C \subset X$, but we remark that the holomorphic tubular neighborhood is generally false: if $S \subset X$ is a compact holomorphic submanifold, there is not necessarily a neighborhood $U \subset X$ of $S$ which is biholomorphic to a neighborhood of the zero section in the total space of the normal bundle $\mathcal{N}_S$. When this is true, $S$ is sometimes said to satisfy the formal neighborhood principle. There is much literature on this subject starting with foundational work of Grauert \cite{Grauert} in the case of codimension one submanifolds. In general codimension, there is work of \cite{ABT09} where the condition for a formal neighborhood involves vanishing of $H^1(S, T_S \otimes {\rm Sym}^k(\mathcal{N}_S^*))$ and $H^1(S, \mathcal{N}_S \otimes {\rm Sym}^{k+1}(\mathcal{N}_S^*))$.

  Returning to $(-1,-1)$ curves $C$ on a Calabi-Yau threefold $X$, it is well-known (e.g. \cite{Kontsevich}) that there exists a neighborhood $U \subset X$ which is biholomorphic to a neighborhood of the zero section in $\mathcal{O}(-1)^{\oplus 2} \rightarrow \mathbb{P}^1$.

  \section{Crossing singularities}

  \subsection{Rolling in the landscape}
Consider the set of all possible Calabi-Yau threefolds. It was noticed by string theorists Green-H\"ubsch and Candelas-Green-H\"ubsch \cite{CGH, GH1, GH2} that there is a way to travel in this landscape by a process known as a conifold transition. This led to the idea of viewing Calabi-Yau threefolds not as isolated objects, but as part of a unified moduli space. We will define a conifold transition momentarily, but let us first state Reid's conjecture. 
  
  \begin{conj}
\cite{Reid} All Calabi-Yau threefolds can be connected by a sequence of conifold transitions.
\end{conj}

Early work on the algebro-geometric foundations of conifold transitions goes back to Clemens \cite{ClemensSolid} and Friedman \cite{Fri86}. A conifold transition, which we will denote by
\[
\hat{X} \rightarrow X_0 \rightsquigarrow X_t,
\]
is a two step process. First, $N$ holomorphic $(-1,-1)$-curves on $\hat{X}$ are contracted to points producing a complex analytic space $X_0$ with $N$ ordinary double point singularities. Second, the complex structure of $X_0$ is deformed in a family such that $X_t$ are smooth complex manifolds for $t \neq 0$.

\begin{rk}
  One can also study an extension of Reid's conjecture where $\hat{X} \rightarrow X_0 \rightsquigarrow X_t$ more generally denotes a birational contraction followed by smoothing. In this generality, $X_0$ admits more complicated singularities than ordinary double points. We focus our attention to conifold transition in this survey, and refer to Gross \cite{Gross} for an alternate version of Reid's conjecture on connecting Calabi-Yau threefolds by general geometric transitions.
\end{rk}

\subsubsection{Topological description}

Broadly speaking, a conifold transition is a type of topological surgery of 6-manifolds. We start with the topological implications and return later to the complex analytic definition. At the level of topology, $N$ sets of the form $(D^4 \times S^2)_i$ are removed from one manifold and replaced by $(S^3 \times D^3)_i$ in the other by gluing along the common boundary $S^3 \times S^2$. Call the first manifold $\hat{X}$ and the second $X_t$. The first Betti number $b_1$ does not change.

Let $[C_i] \in H_2(D^4 \times S^2)_i$ and $[L_i] \in H_3(S^3 \times D^3)_i$ be generators of homology. By the excision principle for Euler characteristic, \cite{ClemensSolid, Rossi, A+B} the topological change is captured by the formula
\begin{equation} \label{top-change}
N = \rho + \mu
\end{equation}
where:
\begin{align*}
N & = \text{number of } (D^4 \times S^2)_i \ \text{removed from} \ \hat{X}, \\
\rho & = \text{number of independent 2-cycles } [C_i] \in H_2(\hat{X},\mathbb{C}), \\
\mu & = \text{number of independent 3-cycles } [L_i] \in H_3( X_t,\mathbb{C}).
\end{align*}
We also have
\begin{align*}
  b_2 &\mapsto b_2 - \rho\\
  b_3 &\mapsto b_3 + 2 \mu.
\end{align*}
We see that as $\hat{X}$ deforms to $X_t$ across a conifold transition $\hat{X} \rightarrow X_0 \rightsquigarrow X_t$, its Betti numbers jump as $b_2 \downarrow$ and $b_3 \uparrow$.

\subsubsection{The definition of a conifold transition}

We now give the definition of a conifold transition following Friedman \cite{Friedman}. Let $\hat{X}$ be a Calabi-Yau threefold. The deformation process of a conifold transition, denoted $\hat{X} \rightarrow X_0 \rightsquigarrow X_t$, is defined as follows:

\begin{enumerate}
\item Find disjoint $(-1,-1)$-curves $C_i \subseteq \hat{X}$.
\item Contract the curves $C_i$ to points $p_i$ to form a complex analytic space $X_0$ with ordinary double point singularities at $p_i$.
\item Realize $X_0$ as the central fiber of a smoothing $\pi: \mathcal{X} \rightarrow \Delta$.
\end{enumerate}

By a smoothing of $X_0$, we mean a proper flat map
\[
\pi: \mathcal{X} \rightarrow \Delta, \quad \pi^{-1}(0) = X_0
\]
with $\mathcal{X}$ a smooth complex fourfold, $\Delta \subseteq \mathbb{C}$ the unit disk, and $\pi^{-1}(t)= X_t$ smooth complex manifolds for $t \neq 0$. Furthermore, we require that a neighborhood $\mathcal{U} \subseteq \mathcal{X}$ of each point $p_i$ is locally analytically isomorphic to the local model
  \[
\pi: \mathcal{V} \rightarrow \mathbb{C}, \quad \pi(z,t) = t,
  \]
  where
  \[
\mathcal{V} = \{  z_1^2+z_2^2+z_3^2+z_4^2 = t \} \subseteq \mathbb{C}^4 \times \mathbb{C}.
  \]
Recall that the notion of a flat map $\varphi: X \rightarrow Y$ generalizes the notion of a holomorphic submersion by allowing singular fibers. A holomorphic map $\varphi: X \rightarrow Y$ between connected complex manifolds is flat if and only if it is an open map. If $\varphi$ is flat in $p$ and the fiber $X_{\varphi(p)}$ is a manifold at $p$, then $\varphi$ is a submersion at $p$. We refer to Fischer's book \cite{Fischer} for these statements and more.
  
  \begin{rk}
  Realizing step (1) is already part of Reid's conjecture. It is an open problem whether every Calabi-Yau threefold $X$ admits $(-1,-1)$-curves, and Reid-Friedman \cite{Reid, Friedman} speculate on whether there is a collection of such curves generating $H_2(X,\mathbb{Z})$.
\end{rk}

We will explain Step (2) in \S \ref{subsec:small}, but we note for now that it follows automatically from Step (1). Step (3) requires a global condition on the configuration of the curves $C_i$. This global condition is Friedman's condition:
\begin{equation} \label{Friedman-cond}
  \sum_{i=1}^N \lambda_i [C_i] =0 \quad {\rm in} \quad H^4(\hat{X},\mathbb{C})
\end{equation}
with each $\lambda_i \neq 0$. That a linear relation such as this must hold can be deduced from the topological change formula \eqref{top-change} which implies $N > \rho$ since $\mu \geq 1$, and with additional work one can show that all $\lambda_i$ are non-zero \cite{Friedman,Kontsevich,RT16}. For extensions of Friedman's condition to nodal Calabi-Yau $n$-folds in higher dimensions $n$, see Rollenske-Thomas \cite{RT16}.

The Friedman-Kawamata-Ran-Tian theorem asserts that if Friedman's condition \eqref{Friedman-cond} is satisfied, then step (3) can be realized. Namely, $X_0$ admits a smoothing. The full statement of the theorem takes the form of a singular version of the BTT theorem for ordinary double point singularities. We state the version here taken from Tian (\cite{Tian} p. 476).

\begin{thm} [Friedman-Kawamata-Ran-Tian Theorem \cite{Fri86, Kawamata, Ran, Tian}] \label{thm:FT} 
  Let $\hat{X}$ be a Calabi-Yau threefold. Let $C_i$ be a collection of $N$ disjoint $(-1,-1)$ curves satisfying Friedman's condition \eqref{Friedman-cond}. Let $\mu: \hat{X} \rightarrow X_0$ be the holomorphic contraction of the $C_i$ with $\mu(C_i)=p_i$ resulting in a complex analytic space $X_0$ with ordinary double point singularities $p_i$.

Then there exists a proper flat map $\pi: \mathcal{X} \rightarrow B$ where $B \subseteq \mathbf{T}$ is a open neighborhood of the origin in $\mathbf{T}=H^1(X_{0,reg},T^{1,0})$. Here $\mathcal{X}$ is a smooth complex manifold and $\pi^{-1}(0)=X_0$.

The vector space of infinitesimal deformations $\mathbf{T}$ can be understood by the exact sequence
  \[
    0 \rightarrow H^1(\hat{X}, T^{1,0}) \rightarrow \mathbf{T} \rightarrow \bigoplus_{i=1}^N \mathbb{C}[p_i] \rightarrow H^4(\hat{X},\mathbb{C})
  \]
  where the last map is $[p_i] \mapsto [C_i]$. Deformations in direction $H^1(\hat{X}, T^{1,0})$ preserve all singular points. Deformations mapped to a vector
  \[
    \sum_i \lambda_i [p_i] \in \bigoplus_{i=1}^N \mathbb{C}[p_i] = \mathbb{C}^N
  \]
  with all $\lambda_i \neq 0$ deform $X_0$ to a family of smooth complex manifolds.
\end{thm}

In summary, one needs to check Friedman's condition \eqref{Friedman-cond} on a configuration of curves on an initial threefold, and then general theory produces a conifold transition to a new threefold. 

\begin{cor} \cite{Friedman, Tian}
Let $\hat{X}$ be a Calabi-Yau threefold, and let $\{ C_i \} \subseteq \hat{X}$ be a collection of disjoint $(-1,-1)$-curves satisfying Friedman's condition \eqref{Friedman-cond}. Then there exists a conifold transition $\hat{X} \rightarrow X_0 \rightsquigarrow X_t$.
\end{cor}

We note that realizing the local model $\mathcal{V} \rightarrow \Delta$ on the global flat family $\mathcal{X} \rightarrow \Delta$ is work of \cite{Kas}.


\subsection{Local model}
Consider the singular point $0 \in V_0$ where
\[
V_0 = \bigg\{ z_1^2+z_2^2+z_3^2+z_4^2 = 0 \bigg\} \subseteq \mathbb{C}^4.
\]
This singularity is sometimes called a conifold singularity, nodal singularity, quadratic singularity, or ordinary double point. It is distinguished by the holomorphic Morse lemma \cite{Atiyah}: a holomorphic function
\[
  f: \mathbb{C}^{n+1} \rightarrow \mathbb{C}, \quad {\rm with} \ f(0)=0, \ Df(0)=0
\]
and non-degenerate holomorphic Hessian matrix $f_{\alpha \beta}(0)$ admits local holomorphic coordinates $w$ near the origin such that
\[
  f(w) \overset{{\rm loc}}{=} \sum_{i=1}^{n+1} w_i^2.
\]
The singularity $0 \in V_0$ can be desingularized in two distinct ways:
  \begin{enumerate}
  \item {\bf By small resolution:} This is a resolution of singularities
    \[
\mu: \hat{V} \rightarrow V_0
    \]
    such that $\hat{V}$ is a smooth complex manifold with $\mu: \hat{V} \backslash E \rightarrow V_0 \backslash \{ 0 \}$ a biholomorphism away from the exceptional set
    \[
      \mu^{-1}(0) = \mathbb{P}^1.
    \]
    The space $\hat{V}$ is biholomorphic to the total space of $\mathcal{O}(-1)^{\oplus 2} \rightarrow \mathbb{P}^1$. We will give further details on this small resolution in \S \ref{subsec:small} below. 

  \item {\bf By smoothing:} The complex analytic space $V_0$ can be deformed into a smooth complex manifold by adding a parameter $t$:
    \[
V_t = \bigg\{ z_1^2+z_2^2+z_3^2+z_4^2 = t \bigg\} \subseteq \mathbb{C}^4.
\]
The spaces $V_t$ for $t \neq 0$ are smooth complex manifolds diffeomorphic to $T^*S^3$. We can insert $V_0$ as the central fiber of a family
\[
\mathcal{V} \overset{\pi}{\rightarrow} \mathbb{C}, \quad \pi(z,t)=t
\]
where
\[
\mathcal{V} = \bigg\{ (z,t) \in \mathbb{C}^4 \times \mathbb{C} : z \in V_t \bigg\}.
\]
We note that $\| z \|^2 \geq |t|$ on $V_t$, and the set $L_t = \{ \| z \|^2 = |t| \} \subseteq V_t$ is called the vanishing cycle. It is diffeomorphic to a 3-sphere, so that $L_t \cong S^3$ and these collapse to a point as $t \rightarrow 0$.
\end{enumerate}

The local model of the conifold transition is then
\begin{equation} \label{localmod}
\hat{V} \overset{\mu}{\rightarrow} V_0 \rightsquigarrow V_t.
\end{equation}
\begin{figure}[htbp]
    \centering
    \includegraphics[width=\textwidth]{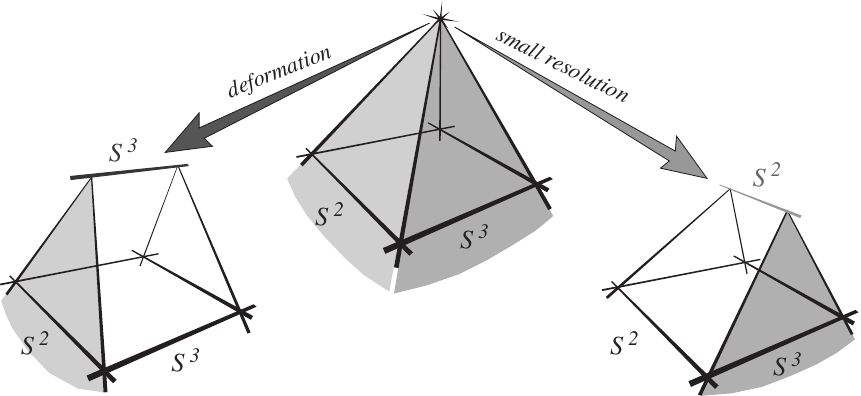}
    \caption{The local model of a conifold transition \cite{CGH,Hubsch}. Illustration taken from H\"ubsch's website [\href{https://tristan.nfshost.com/Research/Conifolds.html}{\ttfamily link}].}
    \label{fig:sample}
\end{figure}
  It should be checked that this description is compatible with the local topological surgery of replacing a neighborhood of the form $D^4 \times S^2$ with one of the form $S^3 \times D^3$. This is done in \cite{ClemensSolid} Lemma 1.11, where in particular the following diffeomorphisms are shown:
  \[
\hat{V} \cong \mathbb{R}^4 \times S^2, \quad V_0 \backslash \{ 0 \} \cong (0,\infty) \times (S^3 \times S^2), \quad V_t \cong S^3 \times \mathbb{R}^3.
  \]
  See also \cite{Rossi} or \cite{A+B} for alternate expositions of these diffeomorphisms. There is a classic diagram depicting the local model of the conifold transition which originates from Figure 1 of Candelas-Green-H\"ubsch \cite{CGH} (see also Figure D.1. of \cite{Hubsch}).

\subsubsection{More on small resolutions} \label{subsec:small}
We now provide more details on the small resolution $\hat{V} \overset{\mu}{\rightarrow} V_0$. We define
\[
\hat{V} = \mathcal{O}(-1)^{\oplus 2} \rightarrow \mathbb{P}^1.
\]
The space $\hat{V} = U \cup \tilde{U}$ is covered by two coordinate charts with coordinates $(u,v,\lambda)$ satisfying the coordinate transformation
\begin{equation} \label{change-coords}
\tilde{\lambda} = \lambda^{-1}, \quad \tilde{u} = \lambda u, \quad \tilde{v} = \lambda v.
\end{equation}
The coordinate $\lambda$ is in the $\mathbb{P}^1$-direction while the $u,v$ are along the fibers of the vector bundle. Let $E \cong \mathbb{P}^1$ be the zero section $\{ u=v=0 \}$. There is a biholomorphism
\[
  \mu: \hat{V} \backslash E \rightarrow V_0 \backslash \{ 0 \},
\]
where
\[
V_0 = \bigg\{ z_1^2+z_2^2+z_3^2+z_4^2 = 0 \bigg\} \subseteq \mathbb{C}^4.
\]
This can be constructed as follows. First, by unitary change of coordinates we identify $V_0$ with
\[
\hat{V}_0 = \{ w_1 w_2 - w_3 w_4 = 0 \} \subseteq \mathbb{C}^4.
\]
Next, we desingularize this space by blowing-up $\mathbb{C}^4$ along $\{w_2=w_4=0 \}$ and taking the proper transform of $\hat{V}_0$. Recall that this means
\[
{\rm Bl} \, \mathbb{C}^4_{\{w_2=w_4=0 \}} = \bigg\{ \left( (w_1,w_2,w_3,w_4), [u,v] \right) \in \mathbb{C}^4 \times \mathbb{P}^1 : (w_2,w_4) \in [u,v] \bigg\}.
\]
The proper transform of $\hat{V}_0$ is biholomorphic to $\hat{V}$. Indeed, over the chart $U = \{ u = 1\}$ there are coordinates $(v,w_2,w_3)$, while over the chart $V = \{ v = 1 \}$ there are coordinates $(u, w_1,w_4)$, and the change of coordinate relation is
\[
 u= v^{-1}, \quad w_1 = v w_3, \quad w_4 = v w_2,
\]
which can be identified with \eqref{change-coords}. The blow-up map of $\hat{V}_0$ induced by $\mu(w,[u]) = w$, is denoted
\[
\mu: \hat{V} \rightarrow \hat{V}_0,
\]
with exceptional set $\mu^{-1}(0) = E \cong \mathbb{P}^1$ so that $\mu: \hat{V} \backslash E \rightarrow \hat{V}_0 \backslash \{ 0 \}$ is a biholomorphism.

\begin{rk}
  The resolution of singularities
  \[
\mu: \hat{V} \rightarrow V_0, \quad \mu^{-1}(0) = \mathbb{P}^1
  \]
  is called a small resolution. It is not unique \cite{Atiyah} as we could have alternatively blown-up $\hat{V}_0 \subseteq \mathbb{C}^4$ along $\{ w_2=w_3=0 \}$.
\end{rk}

\begin{rk}
  If $X_0$ is a complex analytic space with ordinary double point singularities, then a neighborhood of each singularity may be identified with a neighborhood of $0 \in V_0$ and this local procedure defines a small resolution $\mu: \hat{X} \rightarrow X_0$. This construction also defines the contraction of a $(-1,-1)$ curve $C$ on a Calabi-Yau threefold $\hat{X}$ to a singular point of the local form $0 \in V_0$. Indeed, let $C \subseteq \hat{X}$ be a $(-1,-1)$-curve. Then there exists a neighborhood $U \subseteq \hat{X}$ of $C$ biholomorphic to a neighborhood of the zero section in $\hat{V}$, and the local construction defines a contraction $\mu: \hat{X} \rightarrow X_0$.
\end{rk}

\subsection{Examples} Complete intersection Calabi-Yau threefolds can all be connected by conifold transitions; see \cite{GH1, GH2} for work in the string theory literature and \cite{SzWang} in the mathematics literature. From another perspective, rather than connecting known examples, conifold transitions can also be used to construct new examples of projective Calabi-Yau threefolds; see e.g. \cite{BatKr, CandDav}. We list here some simple explicit examples of conifold transitions.

\begin{ex} \label{ex:classic}
  This example can be found in Candelas-Green-H\"ubsch \cite{CGH}. Consider the singular quintic
  \[
X_0 = \bigg\{ Z_3 G(Z_0,\dots,Z_4) + Z_4 H(Z_0,\dots, Z_4) = 0 \bigg\} \subseteq \mathbb{P}^4
  \]
  with $G = Z_3^4 + Z_2^4 - Z_0^4$ and $H=-Z_4^4 - Z_1^4 - Z_0^4$. There are 16 singular points. Near each singularity there exists holomorphic local coordinates such that the local model of the singularity is
  \[
\{ uv -xy = 0 \} \subseteq \mathbb{C}^4,
  \]
  and this is biholomorphic to $V_0 = \{ \sum_i z_i^2 = 0 \} \subseteq \mathbb{C}^4$. We can desingularize $X_0$ in two different ways:
  \begin{itemize}
  \item By small resolution $\hat{X} \rightarrow X_0$. Define $\hat{X} \subseteq \mathbb{P}^4 \times \mathbb{P}^1$ by
    \[
U Z_4 - V Z_3 = 0, \quad U G(Z) + V H(Z) = 0,
    \]
    with $[U,V] \in \mathbb{P}^1$ and $[Z_0, \dots , Z_4] \in \mathbb{P}^4$.
    \item By smoothing
      \[
X_t = \{ Z_3 G + Z_4 H = t Z_0 Z_1 Z_2 Z_3 Z_4 \} \subseteq \mathbb{P}^4.
\]
The $X_t$ are smooth quintics for $t \neq 0$.
    \end{itemize}
    This provides an explicit example of a conifold transition
    \[
\hat{X} \rightarrow X_0 \rightsquigarrow X_t.
    \]
  One can verify that the topological change formula \eqref{top-change} becomes $16=1+15$ by direct calculation of the Betti numbers of $\hat{X}$ and $X_t$, which results in $b_2 \mapsto b_2 -1$ and $b_3 \mapsto b_3 + (2)(15)$. See \cite{Rossi} for a detailed analysis of this example.  \exampleend
  \end{ex}

  \begin{ex} \label{ex:spheres}
    The following example was found by Friedman \cite{Fri86}. Let $\hat{X} \subseteq \mathbb{P}^4$ be a smooth quintic threefold. Then $b_2(\hat{X})=1$, so a pair of disjoint $(-1,-1)$-curves $C_1$, $C_2$ satisfy Friedman's condition \eqref{Friedman-cond}. Choose $C_1$, $C_2$ to generate $H_2(\hat{X},\mathbb{Z})$; see \cite{ClemensIHES, Friedman} for such curves. General theory (Theorem \ref{thm:FT}) gives the existence of a conifold transition
    \[
\hat{X} \rightarrow X_0 \rightsquigarrow X_t.
    \]
We notice that something has gone wrong. The topological change formula \eqref{top-change} implies $b_2(X_t)=0$ and so the complex analytic manifolds $X_t$ cannot admit a K\"ahler structure. In fact, $X_t$ is a simply connected 6-manifold with $H_2(X_t,\mathbb{Z})=0$ and $H_3(X_t,\mathbb{Z})$ torsion-free. By a theorem of Wall \cite{Wall}, the $X_t$ are diffeomorphic to connected sums of $S^3 \times S^3$. \exampleend
\end{ex}

\begin{ex}
  We note that Lu-Tian \cite{LuTian, LuTian2} has constructed examples of conifold transitions from K\"ahler Calabi-Yau threefolds to non-K\"ahler complex manifolds diffeomorphic to $\#_k (S^3 \times S^3)$ for $k \geq 2$. We refer to their work for these explicit examples. \exampleend
  \end{ex}

\begin{rk}
The above examples show that a conifold transition may deform a K\"ahler Calabi-Yau threefold into a non-K\"ahler complex manifold. We will see in the sequel that these complex analytic threefolds inherit some properties of K\"ahler Calabi-Yau manifolds. It is argued by Reid \cite{Reid}, Friedman \cite{Friedman}, Kontsevich \cite{Kontsevich}, Yau \cite{YauNadis,FLY}, Tian \cite{Tian}, etc that these analytic threefolds ought to be included into the mathematical theory developed for the study of projective Calabi-Yau threefolds. Friedman \cite{Fri86} and Reid \cite{Reid} propose a handle-body decomposition for Calabi-Yau threefolds by contraction and smoothing of rational curves spanning $H_2(X,\mathbb{Z})$ so that the resulting threefold is non-K\"ahler and diffeomorphic to a connected sum $\#_k (S^3 \times S^3)$. This program would lead to a sort of uniformization theorem for complex threefolds.
  \end{rk}

\subsubsection{Reversing the arrow}
The directionality of a conifold transition is not standardized in the literature. In all cases, the main feature is the desingularization of ordinary double point singularities by two mechanisms: either small resolution or smoothing of complex structure. Following the conventions established here, we refer to a reverse conifold transition as a deformation where the first step is degeneration of complex structure followed by small resolution:
\[
X_t \rightsquigarrow X_0 \rightarrow \hat{X}.
\]
This direction of travel may also potentially connect a projective Calabi-Yau threefold to a non-K\"ahler threefold. Indeed, it is well-known that Moishezon manifolds \cite{Moi} are not necessarily K\"ahler. 

\begin{ex} Here is a classic example of a quintic with a single ordinary double point. We take the explicit coefficients from \S 1.8 in \cite{CGH}, and this particular example can also be found in H\"ubsch's book (\cite{Hubsch} Ch. D, section D.3.3). Consider the family
  \[
X_t = \bigg\{ Z_5^3 \left( \sum_{i=1}^4 Z_i^2 \right) + \sum_{i=1}^4 a_i Z_i^5 + t Z_5^5 \bigg\} \subseteq \mathbb{P}^4
  \]
  where $a_i$ are non-zero generic constants, and $t \in \mathbb{C}$ is a parameter. As
  \[
X_t \overset{t \rightarrow 0}{\rightsquigarrow} X_0
  \]
  the family of smooth quintics degenerates to a singular space $X_0$ with a single ordinary double point at $P=[0,0,0,0,1]$. Let
  \[
\mu : \hat{X} \rightarrow X_0
  \]
  be a small resolution with $\mu^{-1}(P) = C \cong \mathbb{P}^1$. The formula $N=\rho+\mu$ becomes $1=0+1$, and we see that
  \[
[C] =0 \in H_2(\hat{X},\mathbb{R}).
  \]
  The existence of a K\"ahler metric $\omega$ on $\hat{X}$ would lead to a contradiction:
  \[
0 < \int_C \omega = \int_{\partial \Omega} \omega = 0.
  \]
 We see that the process of degeneration and resolution intertwines K\"ahler and non-K\"ahler manifolds. \exampleend
\end{ex}

\subsection{Summary}
A conifold transition is a mechanism to connect two distinct Calabi-Yau threefolds. We may explore the ways to transfer mathematical structures from one to the other. The following section, Section \S \ref{section:crossing}, will build on this theme from the point of view of differential geometry. For studies on how algebro-geometric or symplectic structures deform through a conifold transition, we refer to e.g. Li-Ruan \cite{LiRuan}, Lee-Lin-Wang \cite{A+B}, Lin-Wang \cite{LinWang}, and references therein.


\section{Geometric structures across singularities} \label{section:crossing}

\subsection{Geometry of local conifold transitions} \label{section:local}
A theme in differential geometry is to understand metric constraints preserved through surgery. We will study conifold transitions in this context. Our starting point is the geometrization of the local model of a conifold transition, and afterwards we will move on to global compact geometries. On the non-compact local model, Ricci-flat metrics were constructed by Candelas-de la Ossa \cite{CdlO90} and we now review their construction. For more details in the presentation style similar to the one taken here, we refer to e.g. \cite{Chuan,CPY,FLY}.

\subsubsection{Small resolution} On the total space
\[
\hat{V} = \mathcal{O}(-1)^{\oplus 2} \rightarrow \mathbb{P}^1,
\]
the Candelas-de la Ossa ansatz \cite{CdlO90} with parameter $a>0$ is
\[
\omega_{co,a} = i \partial \bar{\partial} f_a(r^3) + 4 a^2 \omega_{\rm FS},
\]
where $r$ is a power of the distance function to the zero section $E$ given by
\[
r^3(u,v) = |u|^2_{\omega_{\rm FS}} + |v|^2_{\omega_{\rm FS}},
\]
with $u,v$ fiber coordinates as before. The power of $r$ will be motivated later in \eqref{metric-cone} below. This ansatz is substituted into the Ricci-flat equation where the equation becomes the following ODE
\[
x (f'_a(x))^3 + 6a^2 (f'_a(x))^2 = 1, \quad f_a(x)=f_a(r^3).
  \]
The parameter $a$ measures the volume of $\{r = 0 \} = E \cong \mathbb{P}^1$, so that
\[
{\rm Vol}(\mathbb{P}^1, \omega_{co,a}) \rightarrow 0
\]
as $a \rightarrow 0$. A holomorphic volume form $\hat{\Omega}$ on $\hat{V}$ is given by extending $\mu^* \Omega_0$ across $E$ by Hartog's theorem, where $\mu: \hat{V} \rightarrow V_0$ is a small resolution and $\Omega_0$ is given below in \eqref{V0-holcvol}.

\subsubsection{The singular space} Setting $a=0$, we find the following explicit solution on $\hat{V} \backslash E$:
  \[
\omega_{co,0} = \frac{3}{2} i \partial \bar{\partial} r^2, \quad f_0(x) = \frac{3}{2} x^{2/3}.
  \]
We noted earlier that $\hat{V} \backslash E$ is biholomorphic to the complement of the origin in
\[
V_0 = \bigg\{ \sum_{i=1}^4 z_i^2 = 0 \bigg\} \subseteq \mathbb{C}^4,
  \]
  and the correspondence identifies the radius as $r^3 = \| z \|^2$ on $V_0 \subseteq \mathbb{C}^4$. After rescaling the metric to neglect the factor of 3, the limit of $\omega_{co,a}$ as $a \rightarrow 0$ can be identified with the Ricci-flat geometry
  \[
(V_0, \omega_{co,0}), \quad \omega_{co,0} = \frac{i}{2} \partial \bar{\partial} r^2.
\]
We note that $V_0 \subseteq \mathbb{C}^4$ comes with a scaling action $z \mapsto \lambda z$ and is diffeomorphic to a cone
\[
  V_0 \backslash \{ 0 \} \cong (0,\infty) \times \Sigma, \quad \Sigma = V_0 \cap \{ r = 1 \}
\]
where the link is $\Sigma = S^2 \times S^3$ \cite{CdlO90}. Indeed, if we write $z=x+iy$ then the equation for $V_0$ becomes
\[
\| x \|^2 = \frac{\| z \|^2}{2}, \quad \| y \|^2 = \frac{\| z \|^2}{2}, \quad \langle x,y \rangle =  0.
\]
We interpret the link $\Sigma = \{ \| z \| = 1 \}$ as a sphere bundle $S(T S^3)$ associated to $TS^3 \rightarrow S^3$, and as $TS^3$ is a trivial bundle we have $\Sigma = S^2 \times S^3$.

The metric $\omega_{co,0} = \frac{i}{2} \partial \bar{\partial} r^2$ is a conical metric (see e.g. \cite{Sparks}) as it can be written in polar coordinates
  \begin{equation} \label{metric-cone}
g_{co,0} = dr \otimes dr + r^2 g_\Sigma,
  \end{equation}
where $g_\Sigma$ is a metric on the link $\Sigma$. Thus the Candelas-de la Ossa metric $\omega_{co,0}$ is an example of a Calabi-Yau cone metric. 

\begin{rk}
  The power of $r = \| z \|^{2/3}$ for the Calabi-Yau cone metric can be anticipated by scaling: the K\"ahler-Ricci flat equation is
  \begin{equation}\label{CYeqn}
\omega_{co,0}^3 = i \Omega_0 \wedge \bar{\Omega}_0
  \end{equation}
  where
  \begin{equation} \label{V0-holcvol}
\Omega_0 \overset{{\rm loc}}{=} \frac{dz_1 \wedge dz_2 \wedge dz_3}{z_4}.
  \end{equation}
  The scaling $z \mapsto \lambda z$ leaves \eqref{CYeqn} invariant.
  \end{rk}

  \subsubsection{Smoothing} Next, we equip the smoothings
  \[
V_t = \bigg\{ \sum_{i=1} z_i^2 = t \bigg\} \subseteq \mathbb{C}^4,
  \]
with K\"ahler Ricci-flat metrics. The holomorphic volume form $\Omega_t$ has the same expression as \eqref{V0-holcvol}.  The Candelas-de la Ossa ansatz \cite{CdlO90} for the Calabi-Yau metric is
  \[
\omega_{co,t} = i \partial \bar{\partial} f_t(r^3), \quad r^3 = \| z \|^2,
  \]
  and the Ricci-flat equation on $V_t$ leads to the following ODE for the potential $f_t(x)=f_t(r^3)$:
  \[
(f'_t(x))^3 x + (f'_t(x))^2 f''_t(x) (x^2-|t|^2) = 1/6.
  \]
Setting $t=0$ recovers the solution $f_0(x) = c_0 x^{2/3}$. The submanifold
  \[
L_t = \{ \| z \|^2 = |t| \} \cong S^3
  \]
 is the vanishing cycle $L_t \subseteq V_t$, and it is special Lagrangian with respect to $(V_t,\omega_{co,t},\Omega_t)$. Assuming $t>0$ for simplicity, the special Lagrangian equations take the form
  \[
\omega|_L = 0, \quad {\rm Im} \, \Omega |_L = 0.
\]
The parameter $t$ measures the volume of the special Lagrangian 3-spheres, so that
\[
{\rm Vol}(L_t, \omega_{co,t}) \rightarrow 0
\]
as $t \rightarrow 0$.

  \subsubsection{Summary} The local model of a conifold transition is thus geometrized by K\"ahler Ricci-flat metrics. We have
  \[
(\hat{V}, g_{co,a}) \rightarrow (V_0,g_{co,0}) \leftarrow (V_t,g_{co,t})
  \]
  where convergence is uniform as $a \rightarrow 0$, $t \rightarrow 0$ on compact sets away from the singularities. The convergence is also continuous in the Gromov-Hausdorff sense \cite{FPS24}. The process replaces a holomorphic $\mathbb{P}^1$ by a special Lagrangian $S^3$. For a uniqueness result on K\"ahler Ricci-flat metrics asymptotic to the cone $(V_0,g_{co,0})$, see \cite{ConlonHein}.

\subsection{Global K\"ahler geometry} \label{section:Kahler}
Next, we consider the global setting of a conifold transition of compact Calabi-Yau manifolds. Let $\hat{X}$ be a smooth compact projective Calabi-Yau threefold, and let us degenerate a collection of holomorphic curves to deform $\hat{X}$ by conifold transition
\[
\hat{X} \rightarrow X_0 \rightsquigarrow X_t.
\]
As a first step, we insert the additional hypothesis into our setup that the smoothings $X_t$ happen to be projective. Then both $\hat{X}$ and $X_t$ admit K\"ahler Ricci-flat metrics by Yau's theorem \cite{Yau78}. The question becomes understanding these metrics and their degenerations in families. The study of this problem was instigated by Ruan-Zhang \cite{RuanZhang} and Rong-Zhang \cite{RongZhang}. Let us recall their setup.

\begin{itemize}
\item Let
  \[
\pi: \mathcal{X} \rightarrow \Delta, \quad X_t = \pi^{-1}(t)
  \]
  be a projective smoothing of Calabi-Yau threefolds such that $X_t$ is smooth for $t \neq 0$ and the singular set of $X_0$ consists of finitely many ordinary double points. Let $\mathcal{L} \rightarrow \mathcal{X}$ be an ample line bundle. We consider the family of Calabi-Yau metrics $(X_t,\omega_t)$ defined by
  \[
\omega_t \in c_1(\mathcal{L})|_{X_t}, \quad t \neq 0
  \]
  where $\omega_t$ is K\"ahler Ricci-flat.
\item Let
  \[
\mu: \hat{X} \rightarrow X_0
  \]
  be a small resolution. Let $\hat{\omega}$ be a reference K\"ahler on $\hat{X}$. We consider the family of Calabi-Yau metrics $(\hat{X},\hat{\omega}_a)$ defined by
  \[
\hat{\omega}_a \in a [\hat{\omega}] + \mu^* c_1(\mathcal{L})|_{X_0}, \quad a \neq 0
  \]
  where $\hat{\omega}_a$ is K\"ahler Ricci-flat.
  \end{itemize}

  The theorem of Ruan-Zhang \cite{RuanZhang}, Rong-Zhang \cite{RongZhang}, and J. Song \cite{JSong} is that
  \[
(\hat{X}, \hat{\omega}_a) \rightarrow (X_0, d_0) \leftarrow (X_t,\omega_t)
  \]
  in the Gromov-Hausdorff sense, $(X_0,d_0)$ is a compact metric length space induced by a singular K\"ahler Ricci-flat metric $\omega_0$, and the metrics converge smoothly to $\omega_0$ on compact sets away from the singularities. 

The significance of this result is the manifestation of continuity in a process which is plainly discontinuous topologically as the Betti numbers $b_i$ jump across the conifold transition. Nevertheless, string theory \cite{StromingerBlack,GMS} realizes passing through the conifold singularity in the moduli space as a continuous process. The theorem of \cite{RuanZhang, RongZhang, JSong} is a mathematical counterpart, where the continuity is realized by using families of K\"ahler Ricci-flat metrics and Gromov-Hausdorff convergence.

  The singular Calabi-Yau metric $(X_0,\omega_0)$ exists by Eyssidieux-Guedj-Zeriahi \cite{EGZ}. Better estimate for the singular metric in this context were derived by Hein-Sun \cite{HeinSun}, who showed that the limit $(X_0,g_0)$ inducing the metric space $(X_0,d_0)$ is a singular Calabi-Yau metric with conical singularities, meaning that
  \[
\sup_{ \{ r < \delta \}} |g_0 - c g_{co,0}|_{g_{co,0}} \leq C r^{\lambda}
  \]
  in a small neighborhood of an ordinary double point singularity identified with $V_0 \subseteq \mathbb{C}^4$ and $(V_0,r,\omega_{co,0})$ are as defined above. This is a kind of local uniqueness, where the global metric $\omega_0$ is well-approximated by the explicit Candelas-de la Ossa local model near the singular points. For more results on convergence of global Calabi-Yau metrics to local models, see \cite{ChiuSze, JZhang}.

A consequence of \cite{HeinSun} is that the vanishing 3-cycles can be given the structure of a special Lagrangian 3-sphere with respect to the global Calabi-Yau structure $(X_t,\omega_t,\Omega_t)$. Therefore global Calabi-Yau conifold transitions exchange holomorphic $\mathbb{P}^1$s with a special Lagrangian $S^3$s.

\subsection{Global non-K\"ahler geometry}
We have seen that despite taking initial data $\hat{X}$ to be a projective Calabi-Yau threefold, when deforming across a quadratic singularity
\[
\hat{X} \rightarrow X_0 \rightsquigarrow X_t
\]
the complex threefold $X_t$ will not necessarily remain K\"ahler. We conclude that when considering Calabi-Yau threefolds in families, we are forced to enlarge the space of objects considered by including appropriate limiting manifolds. These limiting objects are certain non-K\"ahler complex manifolds.

We look for geometric structures on $X_t$ to help understand it. The K\"ahler-Ricci flat \eqref{CMA} equation is no longer applicable. Instead, the results of Fu-Li-Yau \cite{FLY} and joint work with Collins and Yau \cite{CPY} find the following structure:

\begin{thm} \label{mainthm} \cite{FLY,CPY} Let $\hat{X}$ be a Calabi-Yau threefold, and let
  \[
    \hat{X} \rightarrow X_0 \rightsquigarrow X_t
  \]
  be a conifold transition. For small enough $t$, the complex manifold $X_t$ admits the structure $(\Omega_t,g_t,h_t)$ solving
\begin{align}
d(|\Omega_t|_{\omega_t} \, \omega_t^2) & = 0 \label{CB} \\
F_{h_t} \wedge \omega_t^2  & = 0 \label{HYM} \\
d \Omega_t & = 0 \label{HOL},
\end{align}
where:
\begin{itemize}
\item $\Omega_t$ is a holomorphic volume form
\item $(g_t,h_t)$ is a pair of hermitian metric on $T^{1,0}X$
\item $\omega_t=  g_t(J_t \cdot, \cdot)$
  \item $F_h= \bar{\partial}(h^{-1} \partial h)$ is the Chern curvature.
\end{itemize}
  Furthermore, near the vanishing cycles $(S^3)_i$, the global metrics $(g_t,h_t)$ satisfy
\begin{align*}
  |g_t - a_i \, g_{co,t}|_{g_{co,t}} &\leq C |t|^\lambda,\\
|h_t - b_i \, g_{co,t}|_{g_{co,t}} &\leq C |t|^\lambda
\end{align*}
for constants $C>1$ and $\lambda>0$ and scaling constants $a_i,b_i$. Here $g_{co,t}$ are the K\"ahler Ricci-flat metrics obtained by Candelas-de la Ossa on the local model $V_t$.
  \end{thm}

We now provide some context for these equations.

\begin{rk}
  On a complex manifold of dimension $n$, a dual notion to the K\"ahler condition $d \omega = 0$ is the balanced metric condition $d \star \omega = 0$, or
  \begin{equation} \label{balanced}
d \omega^{n-1} = 0.
  \end{equation}
  The geometry of such metrics was studied by Michelsohn \cite{Michelsohn}. A theorem of Alessandrini-Bassanelli \cite{AB93} states that the existence of balanced metrics is invariant under modifications. Here $X_t$ is not birational to $\hat{X}$, and Fu-Li-Yau's \cite{FLY} result is the construction of balanced metrics on $X_t$. A conformal change of the metric can be used to go between \eqref{balanced} and \eqref{CB}. 
\end{rk}

\begin{rk}
  The equation
  \[
F_h \wedge \omega^{n-1} = 0
  \]
  is the Hermitian-Yang-Mills equation on the tangent bundle $T^{1,0}X$. This equation was solved for general stable holomorphic vector bundles $E \rightarrow X$ over a K\"ahler manifold by Donaldson-Uhlenbeck-Yau \cite{Donaldson, UY}. On the complex manifolds $\#_k S^3 \times S^3$ from Example \ref{ex:spheres}, \cite{Boz} proved stability of the tangent bundle by algebraic methods.

In principle, this equation could also be solved through conifold transitions for bundles other than the tangent bundle, though it is not evident how a conifold transition creates a stable holomorphic vector bundle $E_t \rightarrow X_t$. There is a proposal for this mechanism in the string theory literature \cite{ABG23}. In the mathematics literature, Chuan \cite{Chuan} solved the Hermitian-Yang-Mills equation with the simplifying assumption that the bundle to carry through the transition is holomorphically trivial in a neighborhood of the contracted curves.
\end{rk}

\begin{rk}
  Though $(X_t,\omega_t,\Omega_t)$ is non-K\"ahler, the vanishing cycles in $X_t$ are also special Lagrangian \cite{CGPY23} in the sense that they are calibrated cycles with respect to a conformal change of metric and the closed 3-form calibration
  \[
{\rm Re} \, e^{-i \hat{\theta}} \Omega
  \]
  for a constant angle $\hat{\theta}$. The concept of such special cycles on (not necessarily K\"ahler) complex manifolds was introduced by Harvey-Lawson \cite{HarveyLawson}. Thus from the perspective of special submanifolds, conifold transitions exchange holomorphic $\mathbb{P}^1$s with special Lagrangian $S^3$s regardless of whether $X_t$ is K\"ahler or not.
\end{rk}

\begin{rk}
Strominger \cite{Strominger} derived equations \eqref{CB}, \eqref{HYM}, \eqref{HOL} from constraints of supersymmetry \cite{BergdeRoo}. A possible interpretation of the theorem is that supersymmetry is preserved through conifold transitions, even though the K\"ahler condition is dropped.
  \end{rk}

\begin{rk}
The system \eqref{CB} \eqref{HYM} \eqref{HOL} is not expected to be the final word on geometrization of conifold transitions. To rigify this structure and for example obtain a finite dimensional moduli space, or alternatively a sort of uniqueness theorem, one needs to impose an additional equation. For example, in heterotic string theory, the heterotic Bianchi identity is coupled to the supersymmetric equations \eqref{CB} \eqref{HYM} \eqref{HOL}, and this gives the finite dimensionality of the moduli space \cite{CMO,CMOS,dlOSvanes,MGFRT-ANN,MGFRT-JDG,PW}. 
\end{rk}

\subsubsection{The geometric system}
We may wonder whether $X_t$ admits Ricci-flat metrics, even though it is non-K\"ahler. From the perspective of string theory, this is not the right equation once the 3-form flux $H$ is non-zero \cite{HullTownsend}. Candelas-Horowitz-Strominger-Witten \cite{CHSW} found the K\"ahler Ricci-flat equation by setting $H=0$, but in general the equations of motion couple the Ricci tensor to a 3-form $H$ and a scalar function $\Phi$ (see e.g. \cite{MartelliSparks} for a modern exposition). 

In the setting of conifold transitions, computing the Ricci tensor of the metric $g$ satisfying \eqref{CB} with $\omega = i g_{\mu \bar{\nu}} \, dz^\mu \wedge d \bar{z}^\nu$ leads to (see e.g. \cite{AMP, PicardSurvey} for this calculation)
    \begin{equation} \label{eom}
R_{mn} + 2 \nabla_m \nabla_n \Phi - {1 \over 4} H_{mpq} H_n{}^{pq} = {1 \over 2} (d H)^p{}_{pmn},
    \end{equation}
    where $H = i (\partial - \bar{\partial})\omega$ and $\Phi = - {1 \over 2} \log |\Omega|_\omega$.

 Although Theorem \ref{mainthm} constructs solutions of \eqref{eom} through conifold transitions, the system of equations of heterotic string theory \cite{Hull,Strominger} require an additional equation on $dH$, the heterotic Bianchi identity, to ensure the cancellation of anomalies \cite{GS}, and this additional equation has not yet been solved through conifold transitions.

  S.-T. Yau has conjectured that the complete set of equations appearing in Strominger's paper \cite{Strominger} is solvable through conifold transitions.

  \begin{conj}
[Yau's conjecture \cite{YauNadis}] Let $\hat{X}$ be a Calabi-Yau threefold, and let
  \[
    \hat{X} \rightarrow X_0 \rightsquigarrow X_t
  \]
  be a conifold transition. For small enough $t$, there exists a triple $(\Omega_t,\omega_t,h_t)$ solving \eqref{CB}, \eqref{HYM}, \eqref{HOL} and
  \begin{equation} \label{AC1}
i \partial \bar{\partial} \omega_t = \alpha'_t ({\rm Tr} \, F_{h_t} \wedge F_{h_t} - {\rm Tr} \, R_{\omega_t} \wedge R_{\omega_t} )
    \end{equation}
    where $R_\omega$ is the Chern curvature of $\omega$ and $\alpha'_t>0$.
    \end{conj}

\begin{rk}
    The conjecture is stated here in its simplest form where the unknown to solve for is a pair of metrics $(g_t,h_t)$ on $T^{1,0} X_t$. Conceivably there could be a setup where $h_t$ is a metric on an auxiliary bundle $E_t \rightarrow X_t$.
  \end{rk}

  \begin{rk}
There are other proposed versions of \eqref{AC1}. In string theory, the equation is a formal expansion about the parameter $\alpha'$ which is only valid in certain regimes, and at higher order in $\alpha'$ the curvature $R_\omega$ should be computed using the Hull connection \cite{Hull, MartelliSparks, MMS}. There is an alternate version of \eqref{AC1} compatible with the formalism of generalized geometry where the equation is
    \begin{equation} \label{AC2}
i \partial \bar{\partial} \omega = \langle F_h \wedge F_h \rangle
    \end{equation}
    where the right-hand side generalizes the special case of $\langle \cdot, \cdot \rangle= {\rm Tr}_{V_0} - {\rm Tr}_{V_1}$ with $F_h$ the curvature of $h=h_0 \oplus h_1$ on $V_0 \oplus V_1$ for a pair of holomorphic vector bundles $V_0$, $V_1$. For the interpretation of this equation as a natural structure in generalized geometry, we refer to \cite{MGF-Survey, MGFM23a, MGFMS, MGFRT-TAMS}.
    \end{rk}

    Assuming this additional equation, either \eqref{AC1} or \eqref{AC2}, can be solved, what are the implications? Returning to the theme of understanding K\"ahler to non-K\"ahler conifold transitions, a major direction for future work in this area is to understand what we learn about $X_t$ from \eqref{AC1}. By associating to $X_t$ the moduli space $\mathcal{M}(X_t)$ of solutions to the equations, we may hope to learn about $X_t$ from $\mathcal{M}(X_t)$. For the implications of these equations in string theory, see e.g. \cite{ADMSE,AIMSSTW,CMO,dlOSvanes,McOSva} and references therein.

\subsubsection{Degenerations of the Hermitian-Yang-Mills equation} We compare this setup to the compact K\"ahler case described in Section \S \ref{section:Kahler}. The analogous non-K\"ahler theorem is best understood as a result on degenerations of Hermitian-Yang-Mills metrics. The background geometry is set by the conformally balanced metrics $\omega$ constructed by Fu-Li-Yau \cite{FLY} along the conifold transition $\hat{X} \rightarrow X_0 \rightsquigarrow X_t$. These satisfy 
\[
{\rm Vol}(C_i, \hat{\omega}_a) \overset{a \rightarrow 0}{\rightarrow} 0, \quad {\rm Vol}(L_t, \omega_{t}) \overset{t \rightarrow 0}{\rightarrow} 0
\]
and the geometry degenerates submanifolds: holomorphic curves $C_i \subseteq \hat{X}$ and special Lagrangian cycles $L_t \subseteq X_t$ are tending to zero volume. The task is to solve the Hermitian-Yang-Mills equation on this degenerating background and analyze its limiting singularities. Combining joint work with T. Collins and S.-T. Yau \cite{CPY}, and B. Friedman and C. Suan \cite{FPS24}, we have the following theorem:

\begin{thm} \cite{CPY, FPS24} Let
  \[
(\hat{X}, \hat{\omega}_a) \rightarrow (X_0,\omega_0) \leftarrow (X_t,\omega_t)
  \]
  be the path of reference Fu-Li-Yau metrics. Then there exists a unique normalized sequence of Hermitian-Yang-Mills metrics $h$ on $T^{1,0}X$ with respect to this degenerating geometry such that
  \[
(\hat{X},\hat{h}_a) \rightarrow (X_0,h_0) \leftarrow (X_t,h_t)
  \]
  in the Gromov-Hausdorff sense. Furthermore the limit $(X_0,\omega_0,h_0)$ is a singular solution to the Hermitian-Yang-Mills equation with conical singularities, in the sense that
  \begin{equation} \label{polydecay}
\sup_{\{ r < \delta \} } |h_0 - c g_{co,0}|_{g_{co,0}} \leq C r^\lambda
  \end{equation}
  near the singularities.
  \end{thm}

  This theorem shows that the Candelas-de la Ossa local model described in \S \ref{section:local} still holds in the global compact case. Namely, the compact metrics $h$ are well-approximated by the K\"ahler Ricci-flat local model $g_{co,0}$ near the singularities, though globally it receives non-K\"ahler corrections. This is an analog of the Hein-Sun theorem \cite{HeinSun} for the Hermitian-Yang-Mills equation in the non-K\"ahler case. The two main steps of proof in \cite{CPY} are:

  \begin{enumerate}
  \item Obtain uniform estimates
    \[
|\hat{h}_a|_{g_{co,a}} + |\hat{h}_a^{-1}|_{g_{co,a}} \leq C
    \]
    near the exceptional curves by adapting the Uhlenbeck-Yau method \cite{UY}. This estimate requires the global condition that $T^{1,0} \hat{X}$ is a stable vector bundle over $\hat{X}$. Taking a limit then yields
    \[
C^{-1} g_{co,0} \leq h_0 \leq C g_{co,0}.
    \]
  \item Upgrade uniform equivalence to polynomial decay \eqref{polydecay}. That is, show that $g_{co,0}^{-1} h_0$ decays to $c {\rm Id}_{T^{1,0}X}$. The toy model for this sort of phenomenon in PDE is the following: suppose $u$ is a scalar function on a cone $V_0$ with $\Delta_{g_{cone}} u =0$ and $|u|\leq C_1$, and show the decay
    \[
|u-u(0)| \leq C_2 r^\lambda.
    \]
    This is elementary for harmonic functions, and that the analogous sort of statement holds for the nonlinear Hermitian-Yang-Mills equation depends on a certain Poincar\'e inequality invoking a stability condition on $T^{1,0} V_0$. We refer to \cite{CPY} for details, and see \cite{JacobWal} for another instance of this technique.
  \end{enumerate}

    \subsubsection{Alternate setups}
   We now note some of the alternate approaches to the geometrization of conifold transitions.

    \begin{itemize}
      \item The equations \eqref{CB}, \eqref{HYM}, \eqref{HOL} have not yet been solved in the reverse direction.
    \[
X_t \rightsquigarrow X_0 \rightarrow \hat{X}
    \]
    Here $X_t$ is a degenerating family of projective threefolds and $\hat{X}$ is a small resolution of singularities. The complex manifold $\hat{X}$, though possibly non-K\"ahler, admits balanced metrics \cite{GiuSpo, GiuSpo2}. 

  \item There are alternate non-K\"ahler equations which may be relevant. One of these is the balanced Chern-Ricci flat equations
    \[
d \omega^2 = 0, \quad |\Omega|_\omega = const.
    \]
    The analysis of these equations was developed in \cite{FWW, Tos15, TW17}. The balanced Chern-Ricci flat equations were solved in \cite{GiuSpo2} across reverse conifold transitions. They are still unsolved in the direction $\hat{X} \rightarrow X_0 \rightsquigarrow X_t$.

  \item There are other options inspired by Type IIB string theory with flux \cite{TsengYau} \cite{Tomasiello}.
     \[
d \omega^2 = 0, \quad i \partial \bar{\partial} \omega = \rho_B.
\]
As noted in \cite{MGF-Survey}, it is not clear how a conifold transition creates a 4-form $\rho_B$ which is Poincar\'e dual to a linear combination of holomorphic curves for a forward conifold transition. In the reverse direction, the small resolution creates holomorphic curves satisfying Friedman's relation and the Type IIB equation may be well-suited. See also \cite{MGFMS} for more on this idea.

\item As $X_t$ cannot be simultaneously complex analytic and symplectic, another alternative is to let go of the complex analytic structure of $X_t$ but preserve the symplectic structure. This point of view was developed by Smith-Thomas-Yau \cite{STY}. In other words, although the initial projective threefold solves $d \omega = 0$ and $d \Omega =0$, across the singularity we may choose to either: preserve $d \omega = 0$ but allow $d \Omega \neq 0$, or preserve $d \Omega = 0$ and allow $d \omega \neq  0$.
  
\end{itemize}

In all cases, the fundamental question is how to use the above geometric structures to constrain the possible manifolds appearing on the other side of a degeneration and resolution.

\subsection{Departure from K\"ahler geometry}
\subsubsection{Complex analytic threefolds} We have seen that conifold transitions can take us out of K\"ahler geometry. Generally speaking, there are sometimes advantages to working with complex analytic threefolds instead of restricting to projective threefolds or manifolds admitting a K\"ahler structure. One such example is J. Pardon’s resolution of the MNOP conjecture \cite{JPardon}. Pardon’s theory of enumeration of holomorphic curves is formulated in the analytic category. His proof of the MNOP conjecture uses the generality of complex analytic families rather than algebraic geometry. The central object is Pardon’s Grothendieck group
\[
  H^*_c(\mathcal{Z},\mathrm{Cpx}_3)
\]
which is the total homology of the double complex
\[
C_{-p}({\rm Cpx}_3, C^q_c(\mathcal{Z})) = \bigoplus_{\mathcal{X} \rightarrow \Delta^p} C_c^q( \mathcal{Z}(\mathcal{X}/\Delta^p))
\]
with direct sum over all complex analytic families $\mathcal{X} \rightarrow \Delta^p$ of threefolds over a $p$-simplex $\Delta^p \subseteq \mathbb{R}^p$. The two differentials are the differential of cohomology with compact support $C_c^*(Z)$ and the other is an alternating sum of restriction to the boundary of the simplex. Here $\mathcal{Z}(\mathcal{X}/\Delta^n)$ is the space of compact 1-cycles lying entirely in the fibers of $\mathcal{X} \rightarrow \Delta^n$.

In Pardon's formalism \cite{JPardon}, curve enumeration theories such as Gromov-Witten theory a la Behrend-Fatechi \cite{Behrend} are homomorphisms out of the Grothendieck group
\[
{\rm GW}:  H^*_c(\mathcal{Z},\mathrm{Cpx}_3) \rightarrow \mathbb{Q}((u))
\]
where powers of $u$ keep track of the genus of the curve. Given a projective threefold $X$ and denoting by $Z(X,\beta)$ the space of curves $C$ with $[C]=\beta$, the constant function ${\bf 1}_{X,Z(X,\beta)}$ defines an element in $[\alpha] \in H^0_c(\mathcal{Z},\mathrm{Cpx}_3)$. This produces an enumerative invariant
\[
{\rm GW}([\alpha]).
\]
Deformation invariance comes from connecting a pair of threefolds by a family $\mathcal{X} \overset{\pi}{\rightarrow} \Delta^1$ to define the same class in $H^0_c(\mathcal{Z},\mathrm{Cpx}_3)$. Indeed, equivalence in Pardon's homology is a way to package a vast generalization of deformation invariance.

This formalism for studying curve enumeration invariants uses cohomology with compact support, and enables the exploration of possibly non-compact spaces of curves on general complex analytic threefolds. In Pardon's application to the MNOP conjecture, the framework allows him to extract open neighborhoods of holomorphic curves and separately deform them to break the curve into a union of isolated rigid curves, and then deduce the general conjecture from the case of local curves \cite{Bryan}.

\subsubsection{The web of threefolds}
Consider the set $\mathcal{W}$ of all complex threefolds connected to projective Calabi-Yau threefolds by conifold transitions. It is unknown what is the precise category of complex threefolds constituting this set. We have seen in Example \ref{ex:spheres} that a conifold transition may connect a K\"ahler Calabi-Yau threefold to a non-K\"ahler complex threefold. It is not necessary to collapse $b_2$ to zero to do this: contract for example only two curves on $\hat{X}$ in Example \ref{ex:classic} to obtain a non-K\"ahler $X_t$ with $b_2(X_t)=1$.

An open question is whether there is a sense in which $\mathcal{W}$ is bounded. One could hope to constrain this space of complex threefolds by equipping them with special geometric structures. Currently, we know that the threefolds $X_t$ with small $t$ linked to a projective Calabi-Yau threefold $\hat{X}$ by conifold transition $\hat{X} \rightarrow X_0 \rightsquigarrow X_t$ satisfy the following properties:

\begin{itemize}
\item There exists a holomorphic volume form $\Omega$ \cite{Friedman}.
  \item There holds $h^{1,0}=h^{0,1}=0$ \cite{Friedman}.
  \item The $\partial \bar{\partial}$-lemma is satisfied \cite{LiChi, Fried-ddbar, TJLee}. Moreover, the Hodge-Riemann bilinear relation holds on $H^{2,1}$.
\item The tangent bundle $T^{1,0}$ is stable with respect to a balanced metric $\omega$ \cite{FLY, CPY}.
\end{itemize}

We see that the analytic threefolds in the web of conifold transitions $\mathcal{W}$ inherit some of the properties of K\"ahler geometry, even though they may or may not actually be K\"ahler.


\begin{rk}
The general theory of balanced Calabi-Yau $\partial \bar{\partial}$-manifolds and their deformation theory has been developed by Wu \cite{Wu06}, Popovici \cite{Popo19}, and Lee \cite{TJLee}.
  \end{rk}

A consequence of Yau's theorem \cite{Yau78} is that K\"ahler Calabi-Yau manifolds have stable tangent bundle. The result of \cite{FLY,CPY} is a generalization of Yau's theorem to the web of threefolds in regions where K\"ahler Ricci-flat metrics cannot exist. As a corollary, we may rule out complex threefolds with unstable $T^{1,0}X$ from appearing as limits of K\"ahler to non-K\"ahler conifold transitions. Recall that given a hermitian metric $\omega$ satisfying $d \omega^2 = 0$ on a complex threefold $X$, stability of $T^{1,0}X$ is the property
  \[
\int_X c_1(F) \wedge \omega^2 < 0
  \]
  for all torsion-free coherent subsheaves $F \subseteq T^{1,0}X$ of rank 1,2. It is possible in this generality that $[\omega^2] = 0$; this is for example the case in Example \ref{ex:spheres} and it was shown in \cite{Boz} that for this particular example $T^{1,0}$ has no holomorphic subbundles.

  It is well-known (e.g. \cite{FriedmanBook}) that stability implies 
  \begin{equation} \label{vanishing}
H^0({\rm End} \, T^{1,0}X)= \{ \lambda \, {\rm Id} : \lambda \in \mathbb{C} \}.
\end{equation}
Stability is also relevant when associating a moduli space to a holomorphic bundle, which can illuminate the topology of the underlying manifold \cite{Donaldson83}. It is hoped that the geometric structures presented in this survey will help improve our understanding of the possible complex analytic threefolds appearing as limits of degenerations and resolutions of Calabi-Yau threefolds.


  \smallskip
  \par {\bf Acknowledgements:} I owe much of my understanding of Calabi-Yau threefolds to frequent discussions with J. Bryan, T. Collins, and S.-T. Yau. I thank the MATRIX-Simons Scholarship program for the opportunity to travel and participate in the MATRIX research program on ``The geometry of moduli spaces in string theory''. Thanks to C. Suan, B. Friedman, R. Friedman, T.-J. Lee, T. H\"ubsch and C.-L. Wang, for discussions and comments.

\end{document}